\documentclass[12pt]{article}
\usepackage{amsfonts}
\newcommand{\Aut}{\mathop{\mathrm{Aut}}\nolimits}
\newcommand{\Sym}{\mathop{\mathrm{Sym}}\nolimits}
\newtheorem{theorem}{Theorem}
\begin{document}

$\phantom{.}$\vskip7cm

\begin{center}
{\Large On the automorphism group of the $m$-coloured random graph}\\[5mm]
{\large Peter J. Cameron and Sam Tarzi}\\[5mm]
School of Mathematical Sciences\\
Queen Mary, University of London\\Mile End Road\\London E1 4NS, UK\\
\texttt{p.j.cameron@qmul.ac.uk}
\end{center}

\begin{abstract}
Let $R_m$ be the (unique) universal homogeneous $m$-edge-coloured countable
complete graph ($m\ge2$), and $G_m$ its group of colour-preserving 
automorphisms. The group $G_m$ was shown to be simple by John Truss. 
We examine the automorphism group of $G_m$, and show that it is the group 
of permutations of $R_m$ which induce permutations on the colours, and 
hence an extension of $G_m$ by the symmetric group of degree $m$. We show 
further that the extension splits if and only if $m$ is odd, and in the
case where $m$ is even and not divisible by~$8$ we find the smallest
supplement for $G_m$ in its automorphism group.
\end{abstract}



\section{Introduction}

Fix an integer $m\ge2$, and 
let $R_m$ be the unique homogeneous universal $m$-edge-colouring of the 
countable complete graph (see Truss~\cite{truss}). (Universality means that
any $m$-edge-coloured finite or countable complete graph is embeddable in
$R_m$, and homogeneity means that every colour-preserving isomorphism between
finite subgraphs extends to an automorphism of $R_m$. The uniqueness is a
special case of Fra\"{\i}ss\'e's theory of countable homogeneous structures.
The graph $R_m$ is the `random $m$-edge-coloured complete graph': that is,
we colour edges independently at random, we obtain $R_m$ with 
probability~$1$. More relevant to us is the fact that the isomorphism class
of $R_m$ is residual in the set of all $m$-coloured complete graphs on a
fixed countable vertex set. See~\cite{oligo} for discussion.)

Let $\Aut(R_m)$ be the group of permutations of the vertex set fixing all 
the colours. Truss~\cite{truss} showed that $\Aut(R_m)$ is a simple group. 

For any permutation $\pi$ of the set of colours, let $R_m^\pi$ be
the graph obtained by applying $\pi$ to the colours. Then $R_m^\pi$ is
universal and homogeneous, and hence isomorphic to $R_m$. This means that,
if $\Aut^*(R_m)$ is the group of permutations of
the vertex set which induce permutations of the colours, then 
$\Aut^*(R_m)$ induces the symmetric group $\Sym(m)$ on the colours; 
so $\Aut^*(R_m)$ is an extension of $\Aut(R_m)$ by $\Sym(m)$. 

The first question we consider here is: when does this extension
split? That is, when is there a complement for $\Aut(R_m)$ in $\Aut^*(R_m)$
(a subgroup of $\Aut^*(R_m)$ isomorphic to $\Sym(m)$ which permutes the
colours)? We also show that $\Aut^*(R_m)$ is the automorphism group of
the simple group $\Aut(R_m)$ (so that the outer automorphism group of
this group is $\Sym(m)$).

\begin{theorem}
The group $\Aut^*(R_m)$ splits over $\Aut(R_m)$ if and only if $m$
is odd.
\end{theorem}

\begin{theorem}
The automorphism group of $\Aut(R_m)$ is $\Aut^*(R_m)$.
\end{theorem}

\section{Proof of Theorem 1}

We show first that the extension does not split if $m$ is
even. Suppose that a complement exists, and let $s$ be an element of this
complement acting as $(1,2)(3,4)\cdots(m-1,m)$ on the colours. Then $s$ maps
the subgraph with colours $1,3,\ldots,m-1$ to its complement. But this is
impossible, since the edge joining points in a $2$-cycle of $s$ has its
colour fixed.

Now suppose that $m$ is odd; we are going to construct a complement.

First, we show that there exists a function $f$ from pairs of distinct
elements of $\Sym(m)$ to $\{1,\ldots,m\}$ satisfying
\begin{itemize}
\item $f(x,y)=f(y,x)$ for all $x\ne y$;
\item $f(xg,yg)=f(x,y)^g$ for all $x\ne y$ and all $g$.
\end{itemize}
To do this, we first define $f(1,y)$ for $y\ne1$ arbitrarily subject to the
condition $f(1,x^{-1}) = f(1,x)^{x^{-1}}$. Note that this condition requires
$f(1,s)^s=f(1,s)$ whenever $s$ is an involution; but this is possible, since
any involution has a fixed point (as $m$ is odd). Then we extend to all pairs
by defining $f(x,y)=f(1,yx^{-1})^x$. A little thought shows that no conflict
arises.

Now we take a countable set of vertices, and let $\Sym(m)$ act semiregularly
on it. Each orbit is naturally identified with $\Sym(m)$; we let $x_i$ denote
the element identified with $x$ in the $i$th orbit, as $i\in\mathbb{N}$
(where orbits are indexed by natural numbers). Then we colour the edges
within each orbit by giving $\{x_i,y_i\}$ the colour $f(x,y)$. For edges
between orbits $i$ and $j$, with $i<j$, we colour $\{x_i,1_j\}$ arbitrarily,
and then give $\{y_i,z_j\}$ the image of the colour of $\{(yz^{-1})_i,1_j\}$
under $z$.

Clearly the group $\Sym(m)$ permutes the colours of the edges
consistently, the same way as it permutes $\{1,\ldots,m\}$.

Next we show that a residual set of the coloured graphs we obtain are
isomorphic to $R_m$. We have to show that, given $m$ finite disjoint sets of
vertices, say $U_1,\ldots, U_m$, the set of graphs containing a vertex $v$
joined by edges of colour $i$ to all vertices in $U_i$  (for $i=1,\ldots,m$)
is open and dense. The openness is clear. To see that it is dense, note that
the $m$ finite sets are contained in the union of a finite number of orbits
(say those with index less than $N$); then, for any $i\ge N$, we are free to
choose the colours of the edges joining these vertices to $1_i$ arbitrarily.

Now by construction, the group $\Sym(m)$ we have constructed meets
$\Aut(R_m)$ in the identity; so it is the required complement.

\bigskip

How close can we get when $m$ is even? The construction in the second part
can easily be modified to show that, if there is a group $G$ which acts as 
$\Sym(m)$ on the set $\{1,\ldots,m\}$, in such a way that all involutions 
in $G$ have fixed points on $\{1,\ldots,m\}$, then $G$ is a supplement for 
$\Aut(R_m)$ in $\Aut^*(R_m)$ (that is, $G.\Aut(R_m)=\Aut^*(R_m)$), and 
$G\cap\Aut(R_m)$ is the kernel of the action of $G$ on $\{1,\ldots,m\}$.
We simply replace $\Sym(m)$ by $G$ in the construction, and in place
of $f(xg,yg)=f(x,y)^g$ we require that $f(xg,yg)=f(x,y)^{g\phi}$, where
$\phi$ is the action of $G$ on $\{1,\ldots,m\}$.

If $m$ is even but not a multiple of $8$, then there is a double cover of 
$\Sym(m)$, for $m$ even, in which the fixed-point-free involutions lift to 
elements of order~$4$. (There are two double covers of $\Sym(n)$ for $n\ge4$,
described in~\cite[Chapter~2]{hh} and called there $\tilde S_m$ and $\hat S_m$.
In $\tilde S_m$, the product of $r$ disjoint transpositions lifts to an 
element of order~$4$ if and only if $r\equiv 1$ or~$2$ mod~$4$, while in 
$\hat S_m$, the condition is that $r\equiv 2$ or~$3$ mod~$4$.)
This
shows that there is a supplement meeting $\Aut(R_m)$ in a group of order~$2$
for $m$ even but not divisible by~$8$.

What happens in the remaining case, when $m$ is a multiple of $8$? Is there
a finite supplement, and what is the smallest such?

\section{Proof of Theorem 2}

Since $\Aut(R_m)$ is primitive and not regular, its
centraliser in the symmetric group is trivial; so $\Aut^*(R_m)$ acts
faithfully on $\Aut(R_m)$ by conjugation. We have to show that there are no
further automorphisms.

A permutation group $G$ of countable degree is said to have the \emph{small
index property} if any subgroup $H$ satisfying $|G:H|<2^{\aleph_0}$ contains
the pointwise stabiliser of a finite set; it has the \emph{strong small index
property} if if any subgroup $H$ satisfying $|G:H|<2^{\aleph_0}$ lies between
the pointwise and setwise stabiliser of a finite set.

\subparagraph{Step 1} $R_m$ has the strong small index property.

This is proved by a simple modification of the  arguments for the case $m=2$.
The small index property is proved by Hodges \emph{et~al.}~\cite{hhls}, 
using a result of Hrushovski~\cite{hrush}; the strong version is a simple
extension due to Cameron~\cite{ssip}.

Hrushovski showed that any finite graph $X$ can be embedded into a finite
graph $Z$ such that all isomorphisms between subgraphs of $X$ extend to
automorphisms of $Z$. Moreover, the graph $Z$ is vertex-, edge- and 
nonedge-transitive. He uses this to construct a generic countable sequence
of automorphisms of $R$. To extend this  to $R_m$ is comparatively
straightforward. It is necessary to work with $(m-1)$-edge-coloured
graphs (regarding the $m$th colour as `transparent'). Now the arguments
of Hodges \emph{et~al.} and Cameron go through essentially unchanged.

\subparagraph{Step 2} Since $\Aut(R_m)$ acts primitively on the vertex set,
with permutation rank $m+1$, the vertex stabilisers are maximal subgroups of
countable index with $m+1$ double cosets. Moreover, any further subgroup of
countable index has more than $m+1$ double cosets.

For let $H$ be a maximal subgroup of countable index. By the strong SIP, $H$
is the stabiliser of a $k$-set $X$. If $g$ maps $X$ to a disjoint $k$-set,
then $HgH$ determines the colours of the edges between $X$ and $X^g$, up to
permutations of these two sets. By universality, there are at least
$m^{k^2}/(k!)^2$ such double cosets. Now it is not hard to prove that
$m^{k^2}/(k!)^2>m$ for $k\ge2$. Hence we must have $k=1$.

\subparagraph{Step 3} It follows that any automorphism permutes the vertex
stabilisers among themselves, so is induced by a
permutation of the vertices which normalises $\Aut(R_m)$. To finish
the proof, we show that the normaliser of $\Aut(R_m)$ in the
symmetric group is $\Aut^*(R_m)$.

This is straightforward. A vertex permutation which normalises
$\Aut(R_m)$ must permute among themselves the $\Aut(R_m)$-orbits on
pairs of vertices, that is, the colour classes; so it belongs to
$\Aut^*(R_m)$.

\end{document}